\documentclass[11pt]{article}

\usepackage{latexsym}
\usepackage{amssymb}

\newtheorem{theorem}{Theorem}

\newtheorem{corollary}{Corollary}

\begin{document}
{
\begin{center}
{\Large\bf
On the truncated multidimensional moment problems in $\mathbb{C}^n$.}
\end{center}
\begin{center}
{\bf Sergey M. Zagorodnyuk}
\end{center}

\noindent
\textbf{Abstract.}
We consider the problem of finding a (non-negative) measure $\mu$ on $\mathfrak{B}(\mathbb{C}^n)$
such that
$\int_{\mathbb{C}^n} \mathbf{z}^{\mathbf{k}} d\mu(\mathbf{z}) = s_{\mathbf{k}}$,  $\forall \mathbf{k}\in\mathcal{K}$.
Here $\mathcal{K}$ is an arbitrary finite subset of $\mathbb{Z}^n_+$, which contains
$(0,...,0)$, and $s_{\mathbf{k}}$ are prescribed complex numbers (we use 
the usual notations for multi-indices).
There are two possible interpretations of this problem. At first, one may consider this problem as an extension of the
truncated multidimensional moment problem on $\mathbb{R}^n$, where the support of the measure $\mu$ is allowed to lie in $\mathbb{C}^n$.
Secondly, the moment problem is a particular case of the truncated moment problem in
$\mathbb{C}^n$, with special truncations. 
We give simple conditions for the solvability of the above moment problem.
As a corollary, we have an integral representation with a non-negative measure for linear functionals
on some linear subspaces of polynomials.

\section{Introduction.}

Throughout the whole paper $n$ means a fixed positive integer.
Let us introduce some notations. As usual, we denote by $\mathbb{R}, \mathbb{C}, \mathbb{N}, \mathbb{Z}, \mathbb{Z}_+$
the sets of real numbers, complex numbers, positive integers, integers and non-negative integers,
respectively.
By $\mathbb{Z}^n_+$ we mean $\mathbb{Z}_+\times \ldots \times\mathbb{Z}_+$, and $\mathbb{R}^n = \mathbb{R}\times \ldots \times\mathbb{R}$,
$\mathbb{C}^n = \mathbb{C}\times \ldots \times\mathbb{C}$,
where the Cartesian products are taken with $n$ copies.
Let $\mathbf{k} = (k_1,\ldots,k_n)\in\mathbb{Z}^n_+$, 
$\mathbf{z} = (z_1,\ldots,z_n)\in\mathbb{C}^n$. Then
$\mathbf{z}^{\mathbf{k}}$ means the monomial $z_1^{k_1}\ldots z_n^{k_n}$, and $|\mathbf{k}| = k_1 + \ldots + k_n$.
By $\mathfrak{B}(\mathbb{C}^n)$ we denote the set of all Borel subsets of $\mathbb{C}^n$.

Let $\mathcal{K}$ be an arbitrary finite subset of $\mathbb{Z}^n_+$, which contains
$\mathbf{0} := (0,...,0)$. 
Let $\mathcal{S} = (s_{\mathbf{k}})_{\mathbf{k}\in\mathcal{K}}$
be an arbitrary set of complex numbers.
We shall consider the problem of finding a (non-negative) measure $\mu$ on $\mathfrak{B}(\mathbb{C}^n)$
such that
\begin{equation}
\label{f1_1}
\int_{\mathbb{C}^n} \mathbf{z}^{\mathbf{k}} d\mu(\mathbf{z}) = s_{\mathbf{k}},\qquad  \forall \mathbf{k}\in\mathcal{K}.
\end{equation}
There are two possible interpretations of this problem. At first, one may consider this problem as an extension of the
truncated multidimensional moment problem on $\mathbb{R}^n$, where the support of the measure $\mu$ is allowed to lie in $\mathbb{C}^n$.
Similar situation is known in the cases of the classical Stieltjes and Hamburger moment problems,
where the support of the measure lies in $[0,+\infty)$ and in $\mathbb{R}$, respectively.
Secondly, and more directly, the moment problem~(\ref{f1_1}) is a particular case of the truncated moment problem in
$\mathbb{C}^n$ (see~\cite[Chapter 7]{cit_985_Curto_Fialkow__Book1}, \cite{cit_983_Kimsey_Woerdeman}, \cite{cit_982_Kimsey_Putinar}), 
with special truncations. These truncations do not include conjugate terms.

It is well known that the multidimensional moment problems are much more complicated than their
one-dimensional prototypes~\cite{cit_1000_Berezansky_1965__Book}, \cite{cit_980_Berg_Christiansen_Ressel__Book}, 
\cite{cit_985_Curto_Fialkow__Book1}, \cite{cit_985_Curto_Fialkow__Book2}, \cite{cit_990_Marshall_Book},
\cite{cit_1000_Schmudgen_Book_2017}. 
An operator-theoretical interpretation of the full 
multidimensional moment problem was given
by Fuglede in~\cite{cit_980_Fuglede}. In general, the ideas of the operator approach to moment problems go back to
the works of Naimark in 1940--1943 and then they were developed by many authors, see historical notes 
in~\cite{cit_15000_Zagorodnyuk_2017_J_Adv_Math_Stud}.
In~\cite{cit_20000_Zagorodnyuk_2019_Concr_Oper} we presented the operator approach to the truncated
multidimensional moment problem in $\mathbb{R}^n$.
Other approaches to truncated moment problems can be found in \cite{cit_985_Curto_Fialkow__Book1}, \cite{cit_985_Curto_Fialkow__Book2},
\cite{cit_1500_Vasilescu_2019}, \cite{cit_15000_Zagorodnyuk_2018_AOT}, \cite{cit_983_Kimsey_Woerdeman},
\cite{cit_982_Kimsey_Putinar}.
Recent results can be also found in~\cite{cit_1510_Yoo_2017}, \cite{cit_981_Idrissi_2019}.

In the case of the moment problem~(\ref{f1_1}) we shall need a modification of the operator approach, since we have no positive
definite kernels here. However, this problem can be passed and we shall come to some commuting bounded operators.
We shall provide a concrete commuting extension for this tuple. Then we apply the dilation theory for commuting contractions
to get the required measure. Consequently and surprisingly, we have very simple conditions for the solvability of the moment problem~(\ref{f1_1})
(Theorem~\ref{t2_1}).
As a corollary, we have an integral representation with a \textit{non-negative} measure for linear functionals $L$
on some linear subspaces of polynomials (Corollary~\ref{c2_1}).

\noindent
{\bf Notations. }
Besides the given above notations we shall use the following conventions.
If H is a Hilbert space then $(\cdot,\cdot)_H$ and $\| \cdot \|_H$ mean
the scalar product and the norm in $H$, respectively.
Indices may be omitted in obvious cases.
For a linear operator $A$ in $H$, we denote by $D(A)$
its  domain, by $R(A)$ its range, and $A^*$ means the adjoint operator
if it exists. If $A$ is invertible then $A^{-1}$ means its
inverse. $\overline{A}$ means the closure of the operator, if the
operator is closable. If $A$ is bounded then $\| A \|$ denotes its
norm.
For a set $M\subseteq H$
we denote by $\overline{M}$ the closure of $M$ in the norm of $H$.
By $\mathop{\rm Lin}\nolimits M$ we mean
the set of all linear combinations of elements from $M$,
and $\mathop{\rm span}\nolimits M:= \overline{ \mathop{\rm Lin}\nolimits M }$.
By $E_H$ we denote the identity operator in $H$, i.e. $E_H x = x$,
$x\in H$. In obvious cases we may omit the index $H$. If $H_1$ is a subspace of $H$, then $P_{H_1} =
P_{H_1}^{H}$ denotes the orthogonal projection of $H$ onto $H_1$.

\section{Truncated moment problems on $\mathbb{C}^n$.}

A solution to the moment problem~(\ref{f1_1}) is given by the following theorem.

\begin{theorem}
\label{t2_1}
Let the moment problem~(\ref{f1_1}) with some prescribed $\mathcal{S} = (s_{\mathbf{k}})_{\mathbf{k}\in\mathcal{K}}$ 
be given. The moment problem~(\ref{f1_1}) has a solution if and only if one of the following
conditions holds:

\begin{itemize}
\item[(a)] $s_{(0,...,0)} > 0$;

\item[(b)] $s_{\mathbf{k}} = 0$, $\forall \mathbf{k}\in\mathcal{K}$. 

\end{itemize}
If one of conditions $(a),(b)$ is satisfied, then there exists a solution $\mu$ with a compact support.
\end{theorem}
\textbf{Proof.}
The necessity part of the theorem is obvious.
Let moment problem~(\ref{f1_1}) be given and one of conditions~$(a)$,$(b)$ holds.
If $(b)$ holds, then $\mu\equiv 0$ is a solution of the moment problem.
Suppose in what follows that $s_{(0,...,0)} > 0$.
Observe that we can include the set $\mathcal{K}$ into the following set:
$$ K_{d} := \{ \mathbf{k}=(k_1,\ldots,k_n)\in\mathbb{Z}^n_+: k_j\leq d,\quad j=1,2,...,n \}, $$
for some large $d\geq 1$. Namely, $d$ may be chosen greater than the maximum value of all possible indices $k_j$ in $\mathcal{K}$.
We now set $s_{\mathbf{k}} := 0$, for $\mathbf{k}\in K_{d}\backslash\mathcal{K}$.
Consider another moment problem of type~(\ref{f1_1}), having a new set of indices $\widetilde{\mathcal{K}} = K_d$.
We are going to construct a solution to this moment problem, which, of course, will be a solution to the original problem.

Consider the usual Hilbert space $l^2$ of square summable complex sequences $\vec c = (c_0,c_1,c_2,...)$,
$\| \vec c \|^2_{l^2} = \sum_{j=0}^\infty |c_j|^2$.
We intend to construct a sequence $\{ x_{\mathbf{k}} \}_{ \mathbf{k}\in \widetilde{\mathcal{K}} }$,
of elements of $l^2$, such that
\begin{equation}
\label{f2_5}
(x_{\mathbf{k}}, x_{\mathbf{0}})_{l^2} = s_{\mathbf{k}},\qquad \mathbf{k}\in \widetilde{\mathcal{K}}.
\end{equation}
The elements of the finite set
$\widetilde{\mathcal{K}}$ can be indexed by a single index, i.e., we assume 
\begin{equation}
\label{f2_9}
\widetilde{\mathcal{K}} = \left\{
\mathbf{k}_0, \mathbf{k}_1,\ldots,\mathbf{k}_\rho 
\right\},
\end{equation}
with $\rho + 1 = |\widetilde{\mathcal{K}}|$, and $\mathbf{k}_0 = (0,...,0)$.
Denote $a := \sqrt{ s_{(0,...,0)} } (> 0)$. 
Set
\begin{equation}
\label{f2_15}
x_{\mathbf{0}} := a \vec e_0,\qquad
x_{\mathbf{k}_j} := \vec e_j + \frac{ s_{\mathbf{k}_j} }{ a } \vec e_0,\qquad j=1,2,...,\rho.
\end{equation}
Here $\vec e_j$ means the vector $\vec c = (c_0,c_1,c_2,...)$ from $l^2$, with $c_j=1$, and $0$'s in other places.
Observe that for this choice of elements $x_{\mathbf{k}}$, conditions~(\ref{f2_5}) hold true.
Moreover, it is important for our future purposes that these elements $x_{\mathbf{k}}$ are linearly
independent.

Consider a finite-dimensional Hilbert space $H := \mathop{\rm Lin}\nolimits \{ x_{\mathbf{k}} \}_{ \mathbf{k}\in \widetilde{\mathcal{K}} }$.
Set
$$ K_{d;l} := \{ \mathbf{k}=(k_1,\ldots,k_n)\in K_d: k_l\leq d-1 \},\qquad l=1,2,...,n. $$
Consider the following operator $W_j$ on $\mathbf{Z}^n_+$:
\begin{equation}
\label{f2_23}
W_j (k_1, \ldots, k_{j-1}, k_j, k_{j+1},\ldots, k_n) = (k_1, \ldots, k_{j-1}, k_j + 1, k_{j+1},\ldots, k_n),
\end{equation}
for $j=1,\ldots,n$. Thus, the operator $W_j$ increases the $j$-th coordinate.
We introduce the following operators $M_j$, $j=1,...,n$, in $H$:
\begin{equation}
\label{f2_25}
M_j \sum_{ \mathbf{k}\in K_{d;j} } \alpha_{\mathbf{k}} x_{\mathbf{k}} = 
\sum_{ \mathbf{k}\in K_{d;j} } \alpha_{\mathbf{k}} x_{ W_j \mathbf{k} },\qquad \alpha_{\mathbf{k}}\in\mathbb{C},
\end{equation}
with $D(M_j) = \mathop{\rm Lin}\nolimits \{ x_{\mathbf{k}} \}_{ \mathbf{k}\in K_{d;j} }$.
Since elements $x_{\mathbf{k}}$ are linearly
independent, we conclude that $M_j$ are well-defined operators.
Operators $M_j$ can be extended to a commuting tuple of bounded operators on $H$.
In fact, consider the following operators $A_j\supseteq M_j$, $j=1,...,n$:
\begin{equation}
\label{f2_29}
A_j \sum_{ \mathbf{k}\in K_{d} } \alpha_{\mathbf{k}} x_{\mathbf{k}} = 
\sum_{ \mathbf{k}\in K_{d;j} } \alpha_{\mathbf{k}} x_{ W_j \mathbf{k} },\qquad \alpha_{\mathbf{k}}\in\mathbb{C}.
\end{equation}
Operators $A_j$ are well defined linear operators on the whole $H$. It can be directly verified that they pairwise commute.
Notice that
\begin{equation}
\label{f2_32}
A_1^{k_1} A_2^{k_2} ... A_n^{k_n} x_{\mathbf{0}} = x_{(k_1,k_2,...,k_n)},\qquad \mathbf{k} = (k_1,...,k_n)\in K_d.
\end{equation}
Relation~(\ref{f2_32}) can be verified using the induction argument.
Since $H$ is finite-dimensional, then
$$ \| A_j \| \leq R,\qquad j=1,2,...,n; $$
for some $R>0$.
Set
\begin{equation}
\label{f2_34}
B_j := \frac{1}{C} A_j,\qquad j=1,...,n,
\end{equation}
where $C$ is an arbitrary number greater than $\sqrt{n} R$.
Then
\begin{equation}
\label{f2_37}
\sum_{j=1}^n \| B_j \|^2 < 1.
\end{equation}
In this case there exists a commuting unitary dilation $\mathcal{U} = (U_1,...,U_n)$ of $(B_1,...,B_n)$,
in a Hilbert space $\widetilde H\supseteq H$, see Proposition 9.2 in~\cite[p. 37]{cit_995_Sz.-Nagy_Book}. Namely, we have:
\begin{equation}
\label{f2_38}
\left.\left( P^{\widetilde H}_H U_1^{k_1} U_2^{k_2} ... U_n^{k_n} \right)\right|_H = B_1^{k_1} B_2^{k_2} ... B_n^{k_n},\qquad
k_1,...,k_n\in\mathbb{Z}_+.
\end{equation}
Moreover, we can choose $\mathcal{U}$ to be minimal, that is, the subspaces $U_1^{k_1} ... U_n^{k_n} H$ will span
the space $\widetilde H$ (see Theorem~9.1 in \cite[p. 36]{cit_995_Sz.-Nagy_Book}):
$$ \widetilde H = \mathop{\rm span}\nolimits \left\{
U_1^{k_1} ... U_n^{k_n} H,\ k_1,...,k_n\in\mathbb{Z}
\right\}. $$
Then the Hilbert space $\widetilde H$ will be separable. 
By~(\ref{f2_34}),(\ref{f2_32}),(\ref{f2_5}),(\ref{f2_38}) we may write for an arbitrary
$\mathbf{k} = (k_1, ...,k_n)\in \widetilde{\mathcal{K}}$:
$$ s_{\mathbf{k}} =
(x_{\mathbf{k}}, x_{\mathbf{0}})_{l^2} = 
(A_1^{k_1} A_2^{k_2} ... A_n^{k_n} x_{\mathbf{0}}, x_{\mathbf{0}})_{l^2} = 
C^{|\mathbf{k}|} (B_1^{k_1} B_2^{k_2} ... B_n^{k_n} x_{\mathbf{0}}, x_{\mathbf{0}})_{l^2} = $$
$$ = C^{|\mathbf{k}|} (U_1^{k_1} U_2^{k_2} ... U_n^{k_n} x_{\mathbf{0}}, x_{\mathbf{0}})_{l^2} = 
( (C U_1)^{k_1} (C U_2)^{k_2} ... (C U_n)^{k_n} x_{\mathbf{0}}, x_{\mathbf{0}})_{l^2} = $$
\begin{equation}
\label{f2_42}
= ( N_1^{k_1} N_2^{k_2} ... N_n^{k_n} x_{\mathbf{0}}, x_{\mathbf{0}})_{l^2},
\end{equation}
where $N_j := C U_j$, $j=1,...,n$.
Applying the spectral theorem for commuting bounded normal operators $N_j$ (or,equivalently, to their commuting
real and imaginary parts), we obtain that
$$ N_j = \int_{\mathbb{C}^n} z_j dF(z_1,...,z_n),\qquad j=1,...,n, $$
where $F(z_1,...,z_n)$ is some spectral measure on $\mathfrak{B}(\mathbb{C}^n)$.
Then
$$ s_{\mathbf{k}} = \int_{\mathbb{C}^n} z_1^{k_1} ... z_n^{k_n} d(F(z_1,...,z_n) x_{\mathbf{0}}, x_{\mathbf{0}})_{l^2},\qquad
\mathbf{k} = (k_1, ...,k_n)\in \widetilde{\mathcal{K}}. $$
This means that $\mu = (F(z_1,...,z_n) x_{\mathbf{0}}, x_{\mathbf{0}})_{l^2}$, is a solution of the moment problem.
Since $N_j$ were bounded, $\mu$ has compact support.
$\Box$

\begin{corollary}
\label{c2_1}
Let $\mathcal{K}$ be an arbitrary finite subset of $\mathbb{Z}^n_+$, which contains
$\mathbf{0}$. Let $L$ be a complex-valued linear functional on 
$$ M = M(\mathcal{K}) := \mathop{\rm Lin}\nolimits \{ z_1^{k_1} ... z_n^{k_n} \}_{ \mathbf{k} = (k_1,...,k_n) \in \mathcal{K} }, $$
such that $L(1)>0$.
Then $L$ has the following integral representation:
\begin{equation}
\label{2_45}
L(p) = \int_{\mathbb{C}^n} p(z_1,...,z_n) d\mu,\qquad \forall p\in M,
\end{equation}
where $\mu$ is a (non-negative) measure $\mu$ on $\mathfrak{B}(\mathbb{C}^n)$, having compact support.
\end{corollary}
\textbf{Proof.} It follows directly from Theorem~\ref{t2_1}. $\Box$

Corollary~\ref{c2_1} can be compared with a well known theorem of Boas, which gives a representation for functionals 
(see~\cite[p. 74]{cit_983_Chihara_Book_1978}). 
It is of interest to consider similar problems with infinite truncations and full moment problems. This will be studied
elsewhere.

}

\noindent
Address:

V. N. Karazin Kharkiv National University \newline\indent
School of Mathematics and Computer Sciences \newline\indent
Department of Higher Mathematics and Informatics \newline\indent
Svobody Square 4, 61022, Kharkiv, Ukraine

Sergey.M.Zagorodnyuk@gmail.com; Sergey.M.Zagorodnyuk@univer.kharkov.ua


\begin{thebibliography}{99}

\bibitem{cit_1000_Berezansky_1965__Book}
Yu. M. Berezansky, Expansions in Eigenfunctions of Selfadjoint Operators, Amer. Math. Soc., Providence, RI, 1968. (Russian edition: Naukova
Dumka, Kiev, 1965).

\bibitem{cit_980_Berg_Christiansen_Ressel__Book}
C.\ Berg, J. P. R. Christensen, P. Ressel, Harmonic Analysis on Semigroups. Springer-Verlag, New York, 1984.

\bibitem{cit_983_Chihara_Book_1978}
T. S. Chihara,  An introduction to orthogonal polynomials. Mathematics and its Applications, 
Vol. 13. Gordon and Breach Science Publishers, New York-London-Paris, 1978. xii+249 pp.

\bibitem{cit_985_Curto_Fialkow__Book1}
R. Curto, L. Fialkow, Solution of the truncated complex moment problem for flat data, 
Memoirs Amer. Math. Soc. 119, no. 568 (1996), x+52 pp.

\bibitem{cit_985_Curto_Fialkow__Book2}
R. Curto, L. Fialkow, Flat extensions of positive moment matrices: Recursively generated relations, 
Memoirs Amer. Math. Soc. 136, no. 648 (1998), x+56 pp.


\bibitem{cit_980_Fuglede}
B.\ Fuglede, The multidimensional moment problem, {\it Expo. Math.},
1 (1983), no.~4, pp.\ 47-65.

\bibitem{cit_981_Idrissi_2019}
K. Idrissi, E. H. Zerouali,  Complex moment problem and recursive relations. Methods Funct. Anal. Topology 25 (2019), no. 1, 15–34.

\bibitem{cit_982_Kimsey_Putinar}
D. P. Kimsey, M. Putinar,
Complex orthogonal polynomials and numerical quadrature via hyponormality.---
{\it Comput. Methods Funct. Theory}, (2018), 1--16.

\bibitem{cit_983_Kimsey_Woerdeman}
D. P. Kimsey,  H. J. Woerdeman,  
The truncated  matrix-valued K-moment problem  on $\mathbb{R}^d$, $\mathbb{C}^d$, and  $\mathbb{T}^d$,
{\it Trans. Am. Math. Soc.} \textbf{365} (10), (2013), 5393--5430.

\bibitem{cit_990_Marshall_Book}
M. Marshall, Positive Polynomials and Sums of Squares, Amer. Math. Soc., Math. Surveys and Monographs, Vol.~146,
2008.

\bibitem{cit_995_Sz.-Nagy_Book}
B. Sz.-Nagy, C. Foias, H. Bercovici, L. K\'erchy, Harmonic analysis of operators on Hilbert space. Second edition. 
Revised and enlarged edition. Universitext. Springer, New York, 2010. xiv+474 pp.

\bibitem{cit_1000_Schmudgen_Book_2017}
K. Schm\"udgen, The moment problem. Graduate Texts in Mathematics, 277. Springer, Cham, 2017. xii+535 pp.

\bibitem{cit_1500_Vasilescu_2019}
F.-H. Vasilescu, Moment problems in hereditary function spaces. Concr. Oper. 6 (2019), no. 1, 64–75.

\bibitem{cit_1510_Yoo_2017}
S. Yoo, Sextic moment problems on 3 parallel lines. Bull. Korean Math. Soc. 54 (2017), no. 1, 299–318.

\bibitem{cit_15000_Zagorodnyuk_2017_J_Adv_Math_Stud}
S. M. Zagorodnyuk, The Nevanlinna-type parametrization for the operator Hamburger moment problem.---
{\it J. Adv. Math. Stud.}, {\bf 10}, No. 2 (2017), 183-199.

\bibitem{cit_15000_Zagorodnyuk_2018_AOT}
S. Zagorodnyuk, On the truncated two-dimensional moment problem.---
{\it Adv. Oper. Theory}, {\bf 3}, no. 2 (2018), 63-74.

\bibitem{cit_20000_Zagorodnyuk_2019_Concr_Oper}
S. M. Zagorodnyuk, The operator approach to the truncated multidimensional moment problem.---{\it Concr. Oper.}, {\bf 6} (2019), no. 1, 1--19.



\end{thebibliography}
\end{document}